\documentclass[12pt]{amsart}
\usepackage{amsfonts}
\usepackage{amssymb}

\usepackage{amsmath}
\usepackage{epsf}
\usepackage{epsfig}
\usepackage{colordvi}
\usepackage{fancybox}

\newtheorem{thm}{Theorem}[section]
\newtheorem{cor}[thm]{Corollary}

\newtheorem{prop}[thm]{Proposition}
\theoremstyle{definition}

\theoremstyle{remark}
\newtheorem{remark}{Remark}[section]

\newcommand{\ds}{\displaystyle}

\begin{document}

\title[Polynomial minimal surfaces of degree five  ]
{Polynomial minimal surfaces of degree five}%

\thanks{2010 {\it Mathematics Subject Classification}: 53A10}
\author{Ognian Kassabov}%
\address{University of Transport, Sofia, Bulgaria}
\email{okassabov@abv.bg}

\keywords{Minimal surface, isothermal parameters, canonical principal parameters, parametric polynomial surface}%

\begin{abstract}
The problem of finding all  minimal surfaces presented in parametric form 
as polynomials is discussed 
by many authors. It is known that the classical Enneper surface is
(up to position in space and homothety)
the only polynomial minimal surface of degree 3 in isothermal parameters.
In higher degrees the problem is quite more complicated. Here we find 
a general form for the functions that generate a polynomial minimal
surface of arbitrary degree via the Weierstrass formula and 
prove that any polynomial minimal surface of degree 5 in isothermal parameters
may be considered as belonging to one of three special families. 
\end{abstract}
\maketitle




\section{Introduction}

The minimal surfaces are a topic of great interest in many areas as 
mathematics, computer science, physics, medicine, architecture.
The reason is that in small areas they have a minimizing property. 

For the applications of minimal surfaces, in particular in computer graphic 
research, it is important to use
minimal surfaces in polynomial form and hence to know all such surfaces
in small degrees. In this direction Cos\'in and Monterde 
\cite{C-M} proved that up to position in space and homothety the
only polynomial minimal surface of degree three in isothermal parameters is the 
classical Enneper surface. The case of degree four is considered in
\cite{K-V}. Polynomial minimal surfaces of degrees five and 
six are studied in \cite{X-W-2} and \cite{X-W-1}, respectively. Theorems
about their coefficients-vectors are found and some examples are considered. Unfortunately 
the systems for the coefficients are very complicated and the general solution is difficult 
to be found. In \cite{X-W-3}, polynomial minimal surfaces of
arbitrary degree constructed on some special functions are studied, and thus some special 
surfaces are proposed. It is remarked that 
in degrees 3 and 5 these surfaces coincide with the Enneper surface and some of the 
surfaces from \cite{X-W-2}, respectively.

In the present paper we first show that a polynomial minimal surface in isothermal
parameters must be generated via the Weierstrass formula with a polynomial and
a rational function (Section 3). Then in Section 4 we determine all polynomial minimal
surfaces of degree five but we do not try to solve the system for the coefficients. 
Instead we use the result from Section 3 and we obtain a list of functions that 
generate via the Weierstrass formula all such surfaces.  The surfaces introduced in \cite{X-W-2} 
belong to one of the obtained families, but these families contain 
many other surfaces as well.

It is natural to ask whether these families contain different surfaces. In general the
problem of comparing surfaces given in different parametric form is very complicated. For 
minimal surfaces such a method is proposed in \cite{OK}. It is based on the
canonical parameters introduced in \cite{GG} and then solving an
ordinary differential equation for finding these parameters.   

When trying to investigate the relation between the families obtained in Section 4  
we cannot use directly the method from \cite{OK}, because we cannot find a simple 
form of the surfaces in canonical parameters. So we change a little the approach 
and we escape solving the differential equation for the transition to canonical 
parameters. As a result we find that the obtained three families contain different
minimal surfaces except in a special case.




\section{Preliminaries}

Let $S$ be a regular surface in the Euclidean space defined by the parametric equation
$$
	\boldsymbol{ x=x}(u,v)=(x_1(u,v),x_2(u,v),x_3(u,v)) , \qquad 
	(u,v) \in U \subset \mathbb{R}^2 .
$$
The derivatives of the vector function $\boldsymbol{x=x}(u,v)$ are usualy denoted by
$	{\boldsymbol x}_u$, $	{\boldsymbol x}_v$, ${\boldsymbol x}_{uu}$, etc. Then the 
coefficients of the first fundamental form are equal to the scalar products
$$
	E={\boldsymbol x}_u^2, \qquad F={\boldsymbol x}_u{\boldsymbol x}_v , \qquad G={\boldsymbol x}_v^2 ,
$$

The unit normal to the surface is the vector field
$$
	{\boldsymbol U} = \frac{{\boldsymbol x}_u\times {\boldsymbol x}_v}{|{\boldsymbol x}_u\times {\boldsymbol x}_v|}
	= \frac{{\boldsymbol x}_u\times {\boldsymbol x}_v}{\sqrt{EG-F^2}} .
$$

In particular, if  $E=G$, $F=0$, then the parameters $(u,v)$ of the surface are called {\it isothermal}.

The coefficients of the second fundamental form of $S$ are given  by 
$$
	L={\boldsymbol U}\,{\boldsymbol x}_{uu} , \qquad 	M={\boldsymbol U}\,{\boldsymbol x}_{uv} , 
	                                          \qquad 	N={\boldsymbol U}\,{\boldsymbol x}_{vv} .
$$

The Gauss curvature $K$ and the mean curvature $H$ of  $S$ are 
defined respectively by
$$
  K=\frac{LN-M^2}{EG-F^2} , \qquad	
  H=\frac{E N-2F M+G L}{2(EG-F^2)} .
$$
Recall that the surface $S$ is called {\it minimal} if its mean 
curvature vanishes identically. In this case it follows easily that
the Gauss curvature is nonpositive. 

The study of minimal surfaces is closely related with some complex curves --
those with isotropic tangent vectors. They are called {\it minimal curves}.
Indeed we have the following construction. 
 
Let $S$ be a minimal surface defined in isothermal parameters.
Then it can be considered as the real part of a minimal curve. More precisely,
let $f(z)$ and $g(z)$ be two holomorphic functions. Define the Weierstrass 
complex curve $\Psi(z)$  by
\begin{equation} \label{eq:2.1}
   {\boldsymbol	\Psi}(z)=\int_{z_0}^z\left( \frac12 f(z)(1-g^2(z)),\frac{i}2 f(z)(1+g^2(z)), f(z)g(z) \right)dz \ . 
\end{equation}
Then $\boldsymbol\Psi(z)$ is a minimal curve and its real
and imaginary parts ${\boldsymbol x}(u,v)$ and ${\boldsymbol y}(u,v)$ are
harmonic functions that  define two  
minimal surfaces in isothermal parametrizations. We say that these two minimal 
surfaces are {\it conjugate}. Moreover, every minimal surface can be 
obtained at least locally as the real (as well as the imaginary) part of  a Weierstrass minimal curve.

For any two conjugate minimal surfaces ${\boldsymbol x}(u,v)$ and ${\boldsymbol y}(u,v)$ 
it is defined the {\it associated family} $\{ S_t \}$, where
$$
	S_t \ \ :\ \qquad \boldsymbol{ x}_t(u,v)={\boldsymbol x}(u,v)\cos t+{\boldsymbol y}(u,v)\sin t \ .
$$
Then for any real number $t$ the surface $S_t$ is also minimal and
has the same first fundamental form as  $S=S_0$.

\vspace{0.2cm}
{\bf Example.} Taking $f(z)=1$, $g(z)=z$, we obtain a Weierstrass minimal curve
whose real part is the classical Enneper surface
$$
	{\boldsymbol x}(u,v)=\left( \frac u2\big(1+v^2-\frac{u^2}3\big),-\frac v2\big(1+u^2-\frac{v^2}3\big),\frac12(u^2-v^2)\right) \ .
$$
It is well known that the Enneper  surface coincide (up to position in space) with 
any surface in its associated family, see e.g. \cite{G-A-S}.

\vspace{0.2cm}

In \cite{GG} Ganchev introduces 
{\it the canonical principal parameters}.  If a surface is  parametrized with 
them, the coefficients of its fundamental forms are given by
\begin{gather*}
	E=\frac1\nu \ , \qquad F=0 \ , \qquad G=\frac1\nu \\
	L=1 \ , \qquad M=0 \ , \qquad N=-1 \ ,
\end{gather*}
where $\nu =\sqrt{-K}$ is the {\it normal curvature} of the surface. 
Actually a surface in canonical principal parametrization is the real part of
a Weierstrass minimal curve generated by some functions $f(z)$, $g(z)$ with 
$f(z)=-1/g'(z)$, i.e. it is the real part of the special Weierstrass curve
$$
	\boldsymbol\Phi(z)=-\int_{z_0}^z  
	\left( \frac{1-g^2(z)}{2g'(z)}, \frac{i(1+g^2(z))}{2g'(z)}, \frac{g(z)}{g'(z)} \right) dz \ .
$$

The canonical principal parameters and the normal curvature play a role similar to that
of the natural parameters and the curvature and torsion of a space curve.  
Namely the following theorem holds:

\begin{thm}\label{T:2.1} {\cite{GG}}  
If a surface is parametrized with canonical principal parameters, 
then its normal curvature $\nu$ satisfies the differential equation
$$
	\Delta \ln \nu+2\nu=0 .
$$
Conversely, for any solution $\nu(u,v)$ of this equation (with $\nu_u\nu_v\ne0$), there exists
an {\bf unique} (up to position in space) minimal surface with normal 
curvature $\nu(u,v)$, $(u,v)$ being canonical principal parameters.
Moreover, the canonical principal parameters $(u,v)$ are determined uniquely up to 
the following transformations 
$$
	\begin{array}{l}
		u=\varepsilon\bar u+a , \\
		v=\varepsilon\bar v+b ,
	\end{array} \qquad \varepsilon=\pm 1 \ , \ a=const., \ b=const. 
$$
\end{thm}

We will also use the following results:

\vspace{0.2cm}
\begin{thm}\label{T:2.2} \cite{OK}
Let the minimal surface $S$ be defined by the
real part of the Weierstrass minimal curve (\ref{eq:2.1}). Any solution 
of the differential equation 
\begin{equation} \label{eq:2.3}
	(z'(w))^2=-\frac1{f(z(w))g'(z(w))}   
\end{equation}
defines a transformation of the isothermal parameters of $S$ to canonical principal 
parameters. Moreover the function $\tilde g(w)$ that defines $S$ via the Ganchev 
formula is given by $\tilde g(w)=g(z(w))$.
\end{thm}

\begin{thm}\label{T:2.3} \cite{OK} Let the holomorphic function $ g(z)$ generate a minimal 
surface in canonical principal parameters, i.e. via the Ganchev formula. 
Then, for an arbitrary complex number $\alpha$, and for an arbitrary real number 
$\varphi$, any of the functions
$$
  e^{i\varphi}\frac{\alpha+g(z)}{ 1-\bar \alpha g(z)}  \ ,  \qquad\qquad 
  \frac{e^{i\varphi}}{g(z)}
$$ 
generates the same surface  in canonical principal parameters (up to position 
in  space). Conversely,  any function that generates this surface (up to position 
in  space) in canonical principal parameters has one of the above forms.
\end{thm}

In sections 4 and 5 we shall consider minimal polynomial surfaces of degree five.
An interesting study of such surfaces is presented in \cite{X-W-2}. First of all it is 
proved that the harmonic condition implies that such a surface must have the form
\begin{equation} \label{eq:2.4}
	\begin{array}{rl}{\boldsymbol r}(u,v)\ = & {\boldsymbol a}\,(u^5-10u^3v^2+5uv^4)+{\boldsymbol b}\,(v^5-10u^2v^3+5u^4v) \\
	                             &+ {\boldsymbol c}\,(u^4-6u^2v^2+v^4)+{\boldsymbol d}\,uv(u^2-v^2) +{\boldsymbol e}\,u(u^2-3v^2)\\
	                             &+{\boldsymbol f}\, v(v^2-3u^2) +{\boldsymbol g}\, (u^2-v^2)+{\boldsymbol h}\, uv 
	                              +{\boldsymbol i}\, u+{\boldsymbol j}\, v+{\boldsymbol k}
	             \end{array}               
\end{equation} 
where $\boldsymbol {a,\, b,\, c,\, d,\, e,\, f,\, g,\, h,\, i,\, j,\, k}$ are coefficient vectors. 
For these coefficients the following holds, see \cite{X-W-2}: 

\begin{thm}\label{T:2.4.} The harmonic polynomial surface  (\ref{eq:2.1})    is minimal 
if and only if its coefficient vectors satisfy the following system of equations
\begin{equation} \label{eq:2.2}
	\left\{\begin{array}{l}
			\boldsymbol{ a}^2={\boldsymbol b}^2   \\
			{\boldsymbol a}.{\boldsymbol b}=0         \\
		  4\,{\boldsymbol a}.{\boldsymbol c}-{\boldsymbol b}.{\boldsymbol d}=0  \\
		  \,{\boldsymbol a}.{\boldsymbol d}+4{\boldsymbol b}.{\boldsymbol c}=0  \\ 
		  16\, {\boldsymbol c}^2-\,{\boldsymbol d}^2+30\,{\boldsymbol a}.{\boldsymbol e}+30\,{\boldsymbol b}.{\boldsymbol f}=0 \\
		  4\, {\boldsymbol d}.{\boldsymbol c} +15\,{\boldsymbol b}.{\boldsymbol e}-15{\boldsymbol a}.{\boldsymbol f}=0          \\
		  9\,{\boldsymbol e}^2-9{\boldsymbol f}^2+16\,{\boldsymbol c}.{\boldsymbol g}-2\,{\boldsymbol d}.{\boldsymbol h}+10\,{\boldsymbol a}.{\boldsymbol i}-10\, {\boldsymbol b}.{\boldsymbol j}=0  \\
		  9\,{\boldsymbol e}.{\boldsymbol f}-4\,{\boldsymbol c}.{\boldsymbol h}-2\,{\boldsymbol d}.{\boldsymbol g}-5\,{\boldsymbol b}.{\boldsymbol i}-5\,{\boldsymbol a}.{\boldsymbol j}=0   \\
		  4\,{\boldsymbol g}^2-\,{\boldsymbol h}^2+6\,{\boldsymbol e}.{\boldsymbol i}+6\,{\boldsymbol f}.{\boldsymbol j}=0   \\
		  2\,{\boldsymbol g}.{\boldsymbol h}-3\,{\boldsymbol f}.{\boldsymbol i}+3\,{\boldsymbol e}.{\boldsymbol j}=0   \\
		  
		  5\,{\boldsymbol a}.{\boldsymbol h}+10\,{\boldsymbol b}.{\boldsymbol g}-12\,{\boldsymbol c}.{\boldsymbol f}+3\,{\boldsymbol d}.{\boldsymbol e}=0   \\
		  5\,{\boldsymbol b}.{\boldsymbol h}-10\,{\boldsymbol a}.{\boldsymbol g}-3\,{\boldsymbol d}.{\boldsymbol f}-12\,{\boldsymbol c}.{\boldsymbol e}=0   \\

		  6\,{\boldsymbol e}.{\boldsymbol g}+3{\boldsymbol f}.{\boldsymbol h}+4{\boldsymbol c}.{\boldsymbol i}-{\boldsymbol d}.{\boldsymbol j}=0    \\
		  6\,{\boldsymbol f}.{\boldsymbol g}-3{\boldsymbol e}.{\boldsymbol h}-{\boldsymbol d}.{\boldsymbol i}-4{\boldsymbol c}.{\boldsymbol j}=0    \\

		  \,{\boldsymbol h}.{\boldsymbol i}+2{\boldsymbol g}.{\boldsymbol j}=0    \\
		  2\,{\boldsymbol g}.{\boldsymbol i}-{\boldsymbol h}.{\boldsymbol j}=0    \\
		  {\boldsymbol i}^2={\boldsymbol j}^2    \\
		  {\boldsymbol i}.{\boldsymbol j}=0  \ .
	\end{array}  \right.
\end{equation}  
\end{thm}

Is seems impossible to find the general solution of the system (\ref{eq:2.2}). So
in \cite{X-W-2} some special solutions are considered  and several interesting 
properties are proved for the obtained surfaces. Using a different approach we shall find 
all polynomial minimal surfaces of degree five.




\setcounter{equation}{0}

\section{Polynomial minimal surfaces of arbitrary degree}

As is said in Introduction, polynomial minimal surfaces of arbitrary degree are studied in 
\cite{X-W-3}. The following construction is proposed. Consider the  functions 
$$
	\begin{array}{l}
		\ds P_n=\sum_{k=0}^{\lceil{\frac{n-1}{2}}\rceil} \, (-1)^k\binom{n}{2k}u^{n-2k}v^{2k}  \\
		\ds Q_n=\sum_{k=0}^{\lfloor{\frac{n-1}{2}\rfloor}} \, (-1)^k\binom{n}{2k+1}u^{n-2k-1}v^{2k+1} \ ,
	\end{array}
$$
where $ \lceil{x}\rceil $ denotes the smallest integer not less than $x$ and $\lfloor{x}\rfloor$ 
denotes the largest integer not greater than $x$.

Then it is proved that for any real number $\omega$ the polynomial surface of degree $n$ defined by
$$
	{\boldsymbol x}(u,v)= \big( -P_n+\omega P_{n-2},\ Q_n+\omega Q_{n-2},\ \frac{2\sqrt{n(n-2)\omega}}{n-1}P_{n-1} \big)
$$
is minimal. Of course this large family is very interesting. But it is important also
to know are these all the possible polynomial minimal surfaces and if not to find other
families. To resolve the last problem we propose the following approach.  

Let 
$$
	S \ :\qquad {\boldsymbol x}(u,v)=(x_1(u,v),x_2(u,v),x_3(u,v)) \ ,
$$
be a polynomial minimal surface of degree $n$ in isothermal parameters. Then $x_i(u,v)$ are polynomials
of degree $\le n$, and at least for one $i=1,2,3$ there is an equality.
Suppose that the parametrization is isothermal and $S$ is defined in an open 
subset of $ \mathbb{R}^2 $, containing $(0,0)$. From Lemma 22.25 in \cite{G-A-S} it
follows that (up to translation) ${\boldsymbol x}(u,v)$ is the real part of the
minimal curve
$$
	{\bf \Psi}(z)=2{\boldsymbol x}\left(\frac z2,\frac z{2i}\right) \ .
$$
So this minimal curve is also polynomial of degree $n$. Then
$$
	{\boldsymbol \Psi'(z)}=(\phi_1(z),\phi_2(z),\phi_3(z))=\left( \frac 12 f(z)(1-g^2(z)),\frac i2 f(z)(1+g^2(z)),f(z)g(z) \right) 
$$
for some functions $f(z)$, $g(z)$ and the 
coordinate functions  $\phi_i(z)$ are polynomials of degree at least 
$\le n-1$ so that at least for one $i$ the degree of $\phi_i(z)$  is exactly $n-1$. 
Hence every one of the functions
$$
	f(z)(1-g^2(z))=2\phi_1(z) \qquad f(z)(1+g^2(z))=-2i\phi_2(z)  \qquad f(z)g(z)=\phi_3(z)
$$
is a polynomial and so
\begin{equation} \label{eq:3.1}
	f(z)=\phi_1(z)-i\phi_2(z)  \qquad  f(z)g^2(z)=-(\phi_1(z)+i\phi_2(z))  \qquad f(z)g(z)=\phi_3(z)  
\end{equation}
are polynomials of degree $\le n-1$ and at least for one $i$ the degree in (\ref{eq:3.1})
is exactly $n-1$. So $f(z)$ is a polynomial of degree  $\le n-1$. Now $(\ref{eq:3.1}) _3$
implies that $g(z)$ is a rational function of the form
\begin{equation} \label{eq:3.2}
	g(z)=\frac{P_p(z)}{Q_q(z)} \ , 
\end{equation}
where  $P_p(z)$ and $Q_q(z)$ are polynomials (of degrees $p$ and $q$, respectively)
with no common zeros. According to $(\ref{eq:3.1}) _2$  the function $f(z)g^2(z)$  
is also a polynomial, so
\begin{equation} \label{eq:3.3}
	f(z)=(Q_q(z))^2R_r(z)  
\end{equation}
with a polynomial $R_r(z)$ (of degree $r$). Moreover since the polynomials
$$
	f(z)= (Q_q(z))^2R_r(z) \hspace{0.5in} f(z)g^2(z)=P^2_p(z)R_r(z)   \hspace{0.5in}  f(z)g(z)=P_p(z)Q_q(z)R_r(z)
$$
are of degree $\le n-1$, then \ $2q+r\le n-1$, $2p+r\le n-1$, $p+q+r\le n-1$ \ and at least once
there is an equality. 

Conversely it is easy to see that for arbitrary polynomials $P_p(z)$, $Q_q(z)$, $R_r(z)$ 
with the above restrictions on $p,q,r$, the functions (\ref{eq:3.2}), (\ref{eq:3.3}) 
generate a minimal polynomial surface of degree $n$ via the Weierstrass formula.  

So we have

\begin{thm}\label{T:3.1} Any polynomial minimal surface of degree $n$ in 
isothermal parameters is
generated via the Weierstrass formula by two functions of the form
$$
	f(z)= (Q_q(z))^2R_r(z) \hspace{1.in} g(z)=\frac{P_p(z)}{Q_q(z)}
$$
where $P_p(z),\, Q_q(z),\, R_r(z)$ are polynomials of degree $p,q,r$, respectively
and  \ $2q+r\le n-1$, $2p+r\le n-1$, $p+q+r\le n-1$ \ with at least one
equality. Conversely any two functions $f(z)$, $g(z)$ with the above form
generate a minimal polynomial surface of degree $n$ via the Weierstrass formula.
\end{thm}




\setcounter{equation}{0}

\section{Consequences for polynomial minimal surfaces of degree five}

With the notations of the previous section we assume $n=5$. Then
$$
	 2q+r\le 4,   \hspace{0.6in}  2p+r\le 4,  \hspace{0.6in}  p+q+r\le 4  
$$ 
with at least one equality. According to the first two equations in Theorem \ref{T:2.4.}
we may assume that (up to a change of the coordinate system) ${\bf a}=(a_1,a_2,0)$,
${\bf b}=(-a_2,a_1,0)$ and hence
\begin{equation} \label{eq:3.6}
	 2q+r\le 4,   \hspace{0.6in}  2p+r\le 4,  \hspace{0.6in}  p+q+r\le 3  \ .
\end{equation} 

Then $(\ref{eq:3.6}) _1$ and $(\ref{eq:3.6}) _2$ imply \ $q\le2$, $p\le 2$, so the
following cases can appear:

1. $p=2$. Then $r=0$.

1.1. $q=0$, i.e. $f(z)=a$, $\displaystyle g(z)=Az^2+Bz+C$, where $a, \, A\ne 0$. 
To obtain another functions defining the surface we can use a consequence of 
the following assertion:

\begin{prop}\label{T:4.1} \cite{OK}
Suppose the pairs $(\tilde f(z),\tilde g(z))$ and $ (f(w),g(w))$ generate 
two minimal surfaces via the Weierstrass formula. Then these surfaces 
coincide (up to translation) iff there exists a function $w=w(z)$, such that
$$
	\tilde  f(z)=f(w(z))w'(z)
  \qquad \text{and} \qquad
	\tilde g(z)=g(w(z)) \ .
$$
\end{prop}

\begin{cor}\label{T:4.2} 
Suppose the pair $ (f(z),g(z))$ generates
a minimal surface via the Weierstrass formula. Then 
for arbitrary numbers $\alpha$ $(\alpha\ne 0)$, $\beta$ the pair 
$$
	\tilde f(z)=\alpha f(\alpha z+\beta) \ ,
  \qquad  \qquad
	 \tilde g(z)=g(\alpha z+\beta)
$$
generates the same minimal surface (up to translation).
\end{cor}

Using Corollary \ref{T:4.2} with  \ $\alpha=\sqrt A\ , \beta=\frac B{2\sqrt A}$ \ 
we can say that the surface is generated by two functions of the form

1.1.  $f(z)=a$, $g(z)=z^2+b$, with $a\ne0$.

Analogously we obtain the cases:

1.2. $p=2$, $q=1$, $r=0$ and \ 
$f(z)=a (z+b)^2$, $\displaystyle g(z)=\frac{cz^2+d}{z+b}$,  with $a,c\ne 0$.

2.1. $p=0$, $q=2$, $ r=0$ and \ 
$f(z)=a(z^2+b)^2$, $\displaystyle g(z)=\frac{1}{z^2+b}$,  with $a\ne 0$;

2.2. $p=1$, $q=2$, $ r=0$ and \ 
$f(z)=a(bz^2+c)^2$, $\displaystyle g(z)=\frac{z+d}{bz^2+c}$, with $a,b\ne 0$;

\vspace{0.3cm}
3. $p=1$, $q=0$, $ r=2$ and \ 
$\displaystyle f(z)=az^2+b$, $g(z)=z+c$, with $a\ne 0$;

\vspace{0.3cm}
4. $p=0$, $q=1$, $ r=2$ and \ 
$ f(z)=(az^2+b)(z+c)^2$, $\displaystyle g(z)=\frac 1{z+c}$, with $a\ne 0$.

We will denote the corresponding surfaces ${\boldsymbol r_{11}}[a,b](u,v)$, 
${\boldsymbol r_{12}}[a,b,c,d](u,v)$ \ etc., respectively.

\begin{remark} The case $p=q=r=1$ is not interesting, because in this case the surface is not of degree 5.
\end{remark}
Now we note that the following can be easily proved:

\begin{prop}\label{T:4.3} 
Consider the surfaces
$$
	S_1\ : \qquad {\boldsymbol x_1}(u,v)={\bf Re}\int_{z_0}^z \left(\frac{f_1(z)}2(1-g_1^2(z)),\frac{if_1(z)}2(1+g_1^2(z)),f_1(z)g_1(z)\right)\,dz \ ,
$$
$$
	S_2\ : \qquad {\boldsymbol x_2}(u,v)={\bf Re}\int_{w_0}^w \left(\frac{f_2(w)}2(1-g_2^2(w)),\frac{if_2(w)}2(1+g_2^2(w)),f_2(w)g_2(w)\right)\,dw \ .
$$
Denote by 
$$
	S_2^s\ : \qquad {\boldsymbol x_2^s}(u,v)={\bf Re}\int_{w_0}^w \left(-\frac{f_2(w)}2(1-g_2^2(w)),\frac{if_2(w)}2(1+g_2^2(w)),f_2(w)g_2(w)\right)\,dw
$$
the surface, symmetric of  $S_2$  about the plane  $Oyz$. Then  $S_1$ and $S_2^s$ coincide if and only if
$$
	f_2(w)=f_1(Z(w))g_1^2(Z(w))Z'(w)  \qquad\qquad g_2(w)=\frac1{g_1(Z(w))}
$$
for some function  $Z(w)$.
\end{prop}

Using this proposition we see that the surfaces from cases 2.1, 2.2 and 4 can be 
viewed as symmetric to those in cases 1.1, 1.2 and 3, respectively.
Consequently  we have

\begin{thm}\label{T:4.4} 
Any polynomial minimal surface of degree 5 in isothermal parameters
coincides up to  position in space and symmetry with a surface generated via the 
Weierstrass formula with the pair of functions

1.1.  $f(z)=a$, $g(z)=z^2+b$, with $a\ne0$.

1.2. 
$f(z)=a (z+b)^2$, $\displaystyle g(z)=\frac{cz^2+d}{z+b}$, with $a,c\ne 0$.

\vspace{0.3cm}
3. 
$\displaystyle f(z)=az^2+b$, $g(z)=z+c$, with $a\ne 0$,

\noindent 
where $a,b,c$ are complex numbers.  
\end{thm}

\begin{remark} 
The family of surfaces introduced in \cite{X-W-2} belongs to the case 1.2. More precisely 
the family from \cite{X-W-2} is defined by
$$
	r(u,v)=\big( X(u,v),Y(u,v),Z(u,v) \big)
$$
with
$$
	\begin{array}{l} 
		X(u,v)=a_1(u^5 - 10 u^3 v^2 + 5 u v^4)-a_2(v^5 - 10 v^3 u^2 + 5 v u^4)        \vspace{1mm} \\
                     \hspace{5cm}   +e_1 u (u^2 - 3 v^2) - e_2 v (v^2 - 3 u^2)    \vspace{1mm} \\
    Y(u,v)=a_2(u^5 - 10 u^3 v^2 + 5 u v^4)+a_1(v^5 - 10 v^3 u^2 + 5 v u^4)        \vspace{1mm} \\ 
                     \hspace{5cm}   + e_2 u (u^2 - 3 v^2) + e_1 v (v^2 - 3 u^2)   \vspace{1mm} \\
    Z(u,v)=\frac{\sqrt{30}}{4} \sqrt{\sqrt{(a_1^2 + a_2^2) (e_1^2 + e_2^2)} - (a_1 e_1 + a_2 e_2)} \,(u^4 - 6 u^2 v^2 + v^4)  \vspace{1mm} \\ 
                    \hspace{2.5cm} -\sqrt{30} \sqrt{ \sqrt{(a_1^2 + a_2^2) (e_1^2 + e_2^2)} + (a_1 e_1 + a_2 e_2)} \, u v (u^2 -v^2)      \ ,           
	\end{array}
$$
where  $a_1,a_2,e_1,e_2$ are real parameters. For $a_2e_1-a_1e_2<0$ the surfaces are minimal.
Such a surface is generated by the Weierstrass formula with
$$
	\begin{array}{l}
		f(z)= 6 (e_1 - i e_2) z^2   \vspace{1mm}  \\
		g(z)= \sqrt{\frac56}\,\frac{\sqrt{ \sqrt{(a_1^2 + a_2^2) (e_1^2 + e_2^2)}-a_1 e_1 - a_2 e_2} +
    i \sqrt{ \sqrt{(a_1^2 + a_2^2) (e_1^2 + e_2^2)}+a_1 e_1 + a_2 e_2}}{e_1 - i e_2} \,z
	\end{array}
$$
so it belongs to the case 1.2 with $b=d=0$.
\end{remark}




\setcounter{equation}{0} 

\section{Relations among the families in Theorem \ref{T:4.4}}

For some special values of the parameters the surfaces from Theorem \ref{T:4.4}
obviously coincide. Namely if $d=-b^2c$ in ${\boldsymbol r_{12}}[a,b,c,d](u,v)$, 
the surface is of type 3. On the other hand even when this equality is not 
satisfied, the corresponding surfaces may look very similar,
as Fig. 5.1 and Fig. 5.2 show.

\begin{tabular}{cc}
		\qquad\qquad\epsfig{file=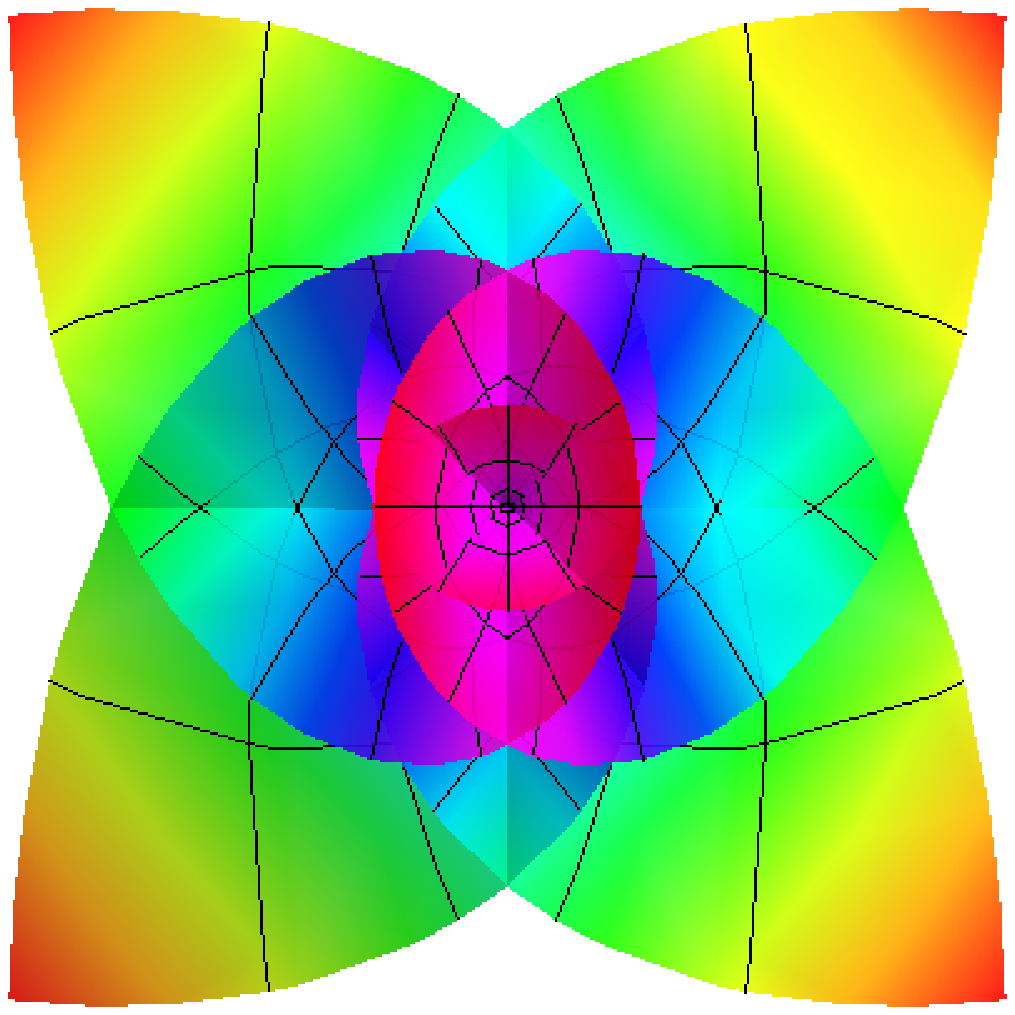,width=0.28\textwidth}   &
    \hspace{2cm}  \qquad \epsfig{file=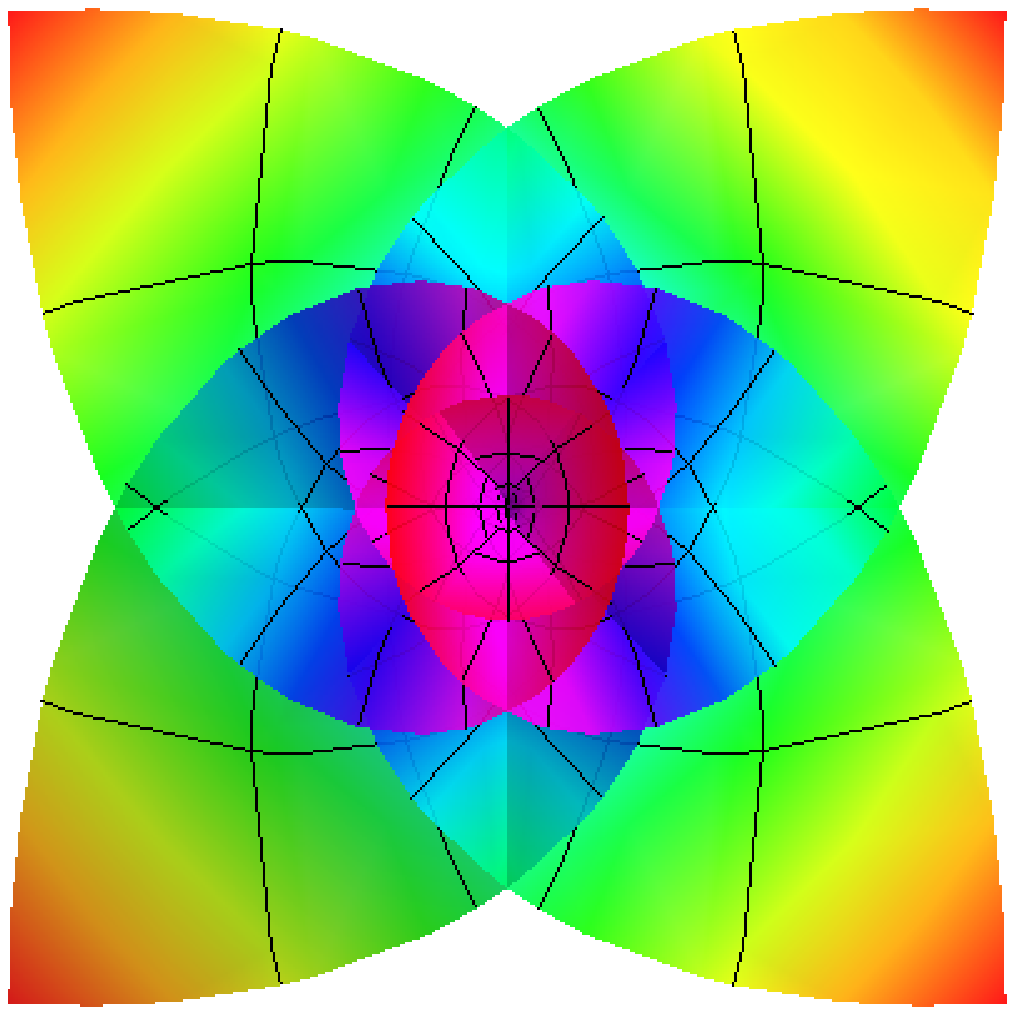,width=0.28\textwidth} \\ 
    \qquad\quad{\bf Fig. 5.1}: ${\boldsymbol r_{12}}[1,0,1,1](u,v)$ & \hspace{2cm}   {\bf Fig. 5.2}: ${\boldsymbol r_3}[1,1,0](u,v)$ \\
\end{tabular} 

We will see that despite the resemblance these two surfaces are different  
as well as that in general the  families  ${\boldsymbol r_{11}}$, ${\boldsymbol r_{12}}$ and
${\boldsymbol r_{3}}$ give different surfaces.

Suppose that a surface ${\boldsymbol r_{12}}[a,b,c,d](u,v)$ generated via 
the Weierstrass formula by the functions
\begin{equation} \label{eq:5.1}
	f_{12}(z)=a(z+b)^2  \hspace{0.8in} g_{12}(z)=\frac{cz^2+d}{z+b}
\end{equation}
coincides (up to position in space) with ${\boldsymbol r_{3}}[A,B,C](u,v)$ generated by
\begin{equation} \label{eq:5.2}
	f_3(z)=Az^2+B  \hspace{0.8in}   g_3(z)=z+C \ .
\end{equation}
Denote $z_{12}(w)$, $z_3(w)$ solutions of the respective equations (\ref{eq:2.3}),
so that (according to Theorem \ref{T:2.3}) the generating functions in canonical principal parameters
\begin{equation} \label{eq:5.3}
  \tilde g_{12}(w)=g_{12}(z_{12}(w))    \hspace{0.8in}  \tilde g_{3}(w)=g_{3}(z_{3}(w)) 
\end{equation}
are related by
\begin{equation} \label{eq:5.4}
	\tilde g_3(w)=e^{i\varphi}\frac{\alpha+\tilde g_{12}(w)}{1-\bar\alpha\, \tilde g_{12}(w)} 
	\qquad {\rm or} \qquad \tilde g_{3}(w)=\frac{e^{i\varphi}}{ \tilde g_{12}(w)} \ .
\end{equation}
We will consider only the first possibility. The second can be considered analogously. 
Note that according to the equation (\ref{eq:2.3}) the functions $z_{12}(w)$ and $z_3(w)$
are related by
\begin{equation} \label{eq:5.5}
	f_{12}(z_{12}(w))g'_{12}(z_{12}(w))(z'_{12})^2=f_{3}(z_{3}(w))g'_{3}(z_{3}(w))(z'_{3})^2 \ .
\end{equation}
From the last equality, using (\ref{eq:5.1})--(\ref{eq:5.4}) and comparing the coefficients of
$z_{12}(w)$ (note that $z_{12}(w)$ may not be constant) we may derive 
$$
	\alpha=0 \qquad a=A c^3 e^{4i\varphi} \qquad C + 2\, b\, c\, e^{i\varphi}=0 \qquad B=0 \qquad b^2 c + d=0 .
$$

So the surfaces    can coincide only if $ b^2 c + d=0 $. In particular
the surfaces defined by ${\boldsymbol r_{12}}[1,0,1,1](u,v)$ and
${\boldsymbol r_3}[1,1,0](u,v)$ are different despite the resemblance in
Figures 5.1 and 5.2. Actually in a smaller neighborhood of $(u,v)=(0,0)$ 
(with the same viewpoint as for Figures 5.1 and 5.2) the difference 
is clear, see Figures 5.3 and 5.4. 

\begin{tabular}{cc}
		\epsfig{file=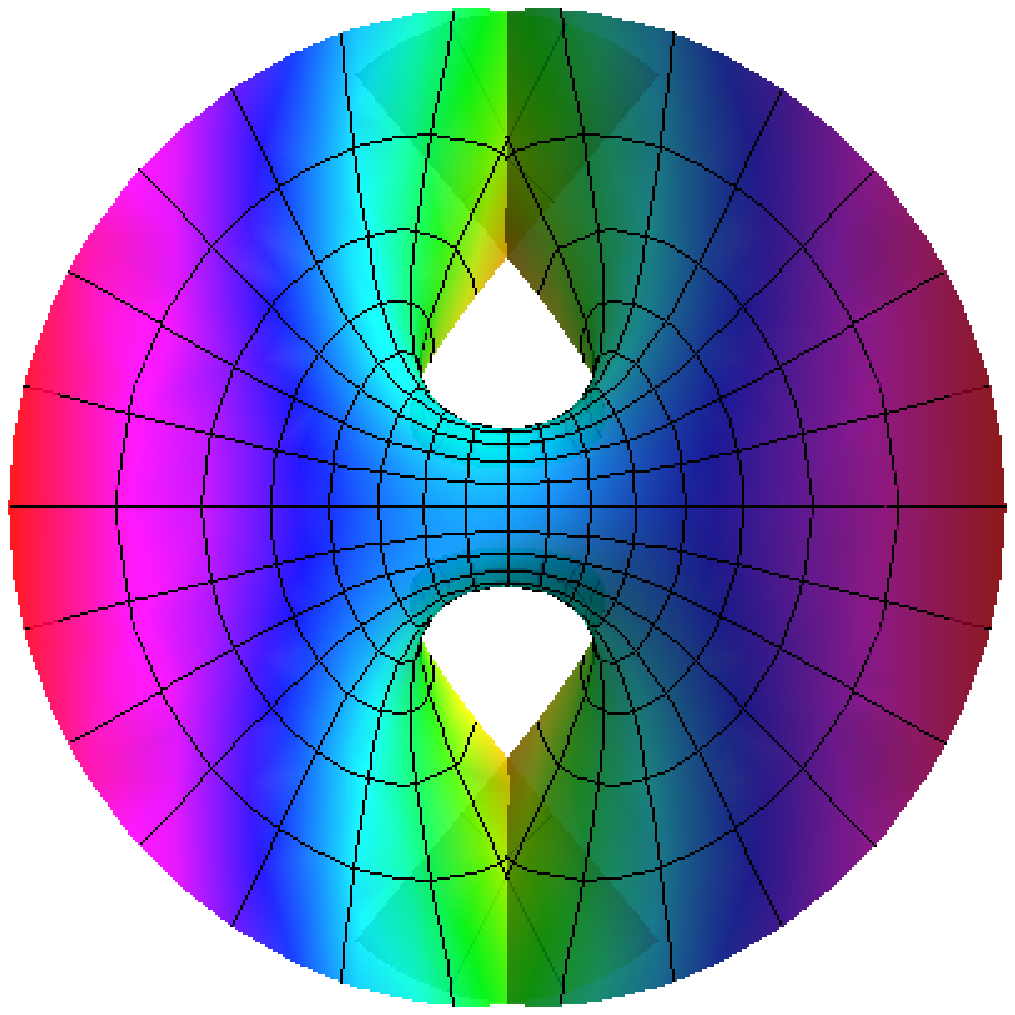,width=0.28\textwidth}   &
      \qquad \epsfig{file=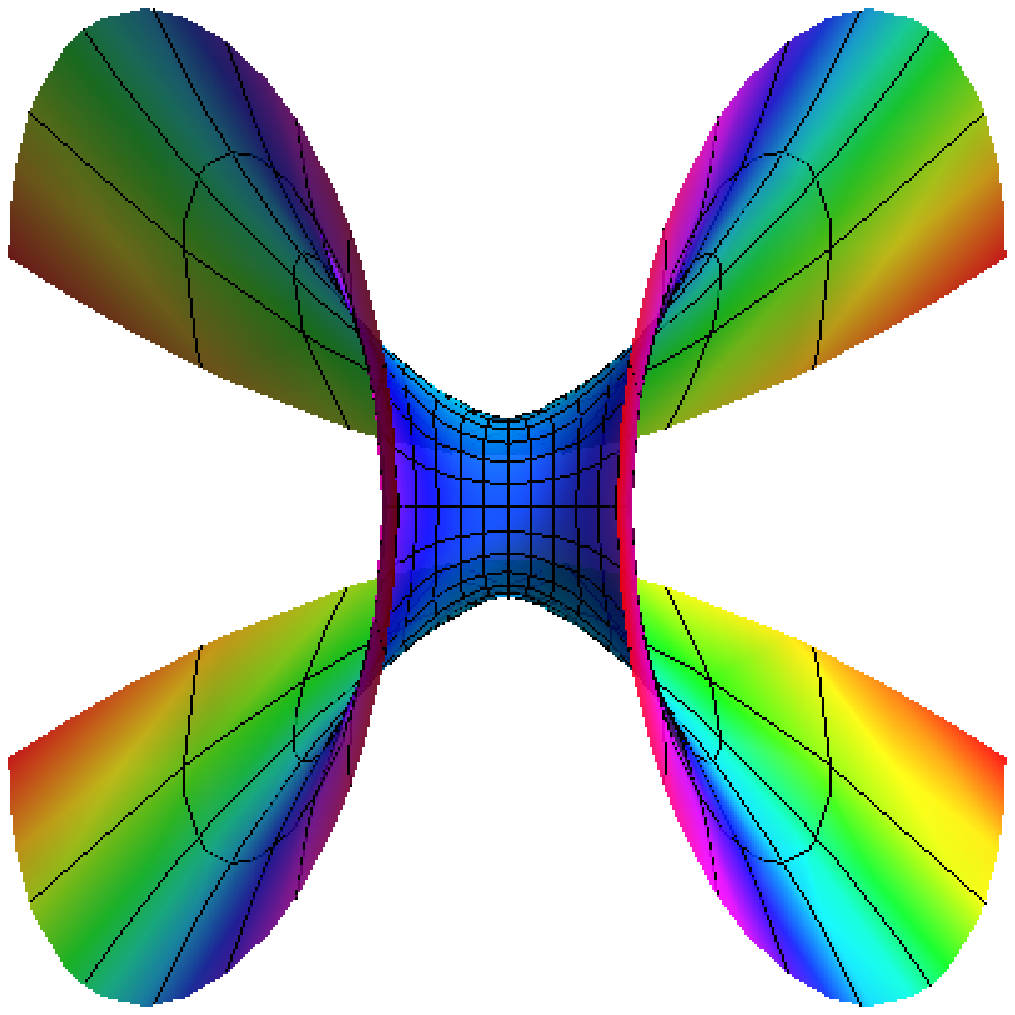,width=0.28\textwidth} \\ 
   \hspace{-0.3cm} {\bf Fig. 5.3}: ${\boldsymbol r_{12}}[1,0,1,1](u,v),\ |u|,|v|\le 1$ 
          & \hspace{0.6cm}   {\bf Fig. 5.4}: ${\boldsymbol r_3}[1,1,0](u,v),\ |u|,|v|\le 1.2$ \\
\end{tabular} 

\vspace{.2cm}We use a similar idea to investigate a possible 
coincidence of surfaces from  cases 1.1 and 1.2.
Let the surface $S_{11} \, :\, {\boldsymbol r_{11}}[A,B](u,v)$ be generated via
the Weierstrass formula by
\begin{equation} \label{eq:5.6}
	f_{11}(z)=A  \hspace{0.8in}   g_{11}(z) = z^2 + B\ .
\end{equation}
and suppose that it coincides with $S_{12} \ :\ {\boldsymbol r_{12}}[a,b,c,d](u,v)$.
Then some functions $\tilde g_{11}$ and $\tilde g_{12}$ that generate them 
in canonical parameters are related by
\begin{equation} \label{eq:5.7}
	\tilde g_{11}(w)=e^{i\varphi}\frac{\alpha+\tilde g_{12}(w)}{1-\bar\alpha\, \tilde g_{12}(w)} 	
	\qquad {\rm or} \qquad \tilde g_{11}(w)=\frac{e^{i\varphi}}{ \tilde g_{12}(w)} \ .
\end{equation}
As before, we consider only the first possibility.
Denote  $z_{11}(w)$, $z_{12}(w)$ respective solutions of the equation (\ref{eq:2.3}).
Then 
\begin{equation} \label{eq:5.8}
	\tilde g_{11}(w)=(z_{11}(w))^2+B
\end{equation}
and hence
\begin{equation} \label{eq:5.9}
	2z_{11}(w)z'_{11}=\tilde g'_{11}(w) \ .
\end{equation}
On the other hand analogously to (\ref{eq:5.5}) 
$$
	\big(f_{11}(z_{11}(w))g'_{11}(z_{11}(w))\big)^2(z'_{11})^4=\big(f_{12}(z_{12}(w))g'_{12}(z_{12}(w))\big)^2(z'_{12})^4 
$$
holds.
Applying (\ref{eq:5.6})--(\ref{eq:5.9}) we can find the left hand side of the 
above equality as a function of $z_{12}(w)$. Then looking at the coefficients of 
$z_{12}(w)$ we conclude that this equality implies a contradiction. 
So a surface $S_{11}$ can not coincide with a surface $S_{12}$. 

Applying the same arguments we may prove that a surface $S_{11}$ can not
coincide with a surface $S_3$. So we have
 
\begin{thm}\label{T:5.1} 
The families from Theorem \ref{T:4.4} contain different
surfaces except if \ $b^2c+d=0$ in case 2.1 and then the surface belongs also
to the case 3.  
\end{thm}

\begin{remark} 
 The surfaces generated via the 
Weierstrass formula by the pairs of functions
$$
	(f(z),g(z))   \qquad {\rm and } \qquad (Cf(z),g(z))
$$
are homothetic for any positive real number $C$. On the other hand, if $C\ne 0$ 
is not real, these surfaces are different in general. 
More precisely let $C=|C|e^{i\varphi}$ for a real number $\varphi$.
The pairs $(f(z),g(z))$ and $(|C|f(z),g(z))$ generate two homothetic surfaces, but
the surface generated by  $(e^{i\varphi}|C|f(z),g(z))$ belongs to the associated
family of the surface generated by $(|C|f(z),g(z))$. Thus we see that if a surface 
belongs to a family from Theorem \ref{T:4.4} then its homothetic surfaces, as well
as their associated surfaces, belong to the same family.
\end{remark}

\end{document}